\documentclass[12pt]{article}
\usepackage{latexsym,amsfonts,amssymb,mathrsfs}
\usepackage[dvips]{color}
\usepackage{colordvi,multicol}

\def \cal{\mathcal}

\newtheorem{thm}{Theorem}[section]
\newtheorem{cor}[thm]{Corollary}

\newtheorem{pro}[thm]{Proposition}
\newtheorem{defi}[thm]{Definition}
\newtheorem{rem}[thm]{Remark}
\newtheorem{exa}[thm]{Example}

\begin{document}
\title{\bf Convergence rates in the law of large numbers and new  kinds of convergence of  random variables}
\author{Ze-Chun Hu$^a$  and  Wei Sun$^b$\\ \\
 {\small $^a$ College of Mathematics, Sichuan  University,  China}\\
 {\small $^b$ Department of
Mathematics and Statistics, Concordia University, Canada}}

\maketitle
\date{}

% \centerline{Ze-Chun Hu\thanks{Corresponding author.}} \centerline{\small College
%of Mathematics, Sichuan University, Chengdu, 610064, China}
%\centerline{\small E-mail: zchu@scu.edu.cn}
%\vskip 0.7cm \centerline{Qian-Qian Zhou} \centerline{\small Department of
%Mathematics, Nanjing University, Nanjing, 210093, China} \centerline{\small
%E-mail: 15266479708@163.com}

%\vskip 1.4cm

\vskip 0.5cm \noindent{\bf Abstract}\quad In this paper, we first
study convergence rates in the law of large numbers for
independent and identically distributed random variables. We obtain a strong $L^p$-convergence version and a
strongly almost sure convergence version of the law of large
numbers. Second, we investigate several new  kinds of convergence of random variables
and discuss their relations and properties.

\smallskip


\smallskip


\tableofcontents

\section{Introduction}

It is well known that limit theorems play
 an important role in the probability theory and statistics. Let
$(\Omega,\mathcal{F},P)$ be a probability space and $\{X,X_n,n\geq
1\}$ be a sequence of random variables. We have different kinds of
convergence:

\begin{itemize}

\item $\{X_n,n\geq 1\}$  is said to almost surely converge to $X$,
if there exists a set $N\in \mathcal{F}$ such that $P(N)=0$ and
$\forall\omega\in \Omega\backslash N,\
\lim_{n\to\infty}X_n(\omega)=X(\omega)$, which is denoted by
$X_n\stackrel{a.s.}{\longrightarrow} X$ or $X_n\to X\ a.s.$.

\item $\{X_n,n\geq 1\}$  is said to converge to $X$ in probability, if for any $\varepsilon>0$, $\lim_{n\to\infty}P(\{|X_n-X|\geq \varepsilon\})=0$,  which is denoted by $X_n\stackrel{P}{\longrightarrow} X$.

\item $\{X_n,n\geq 1\}$  is said to $L^p$-converge to $X$ $(p>0)$
if $\lim_{n\to\infty}E[|X_n-X|^p]=0$, which is denoted by
$X_n\stackrel{L^p}{\longrightarrow} X$.

\item $\{X_n,n\geq 1\}$  is said to $L^{\infty}$-converge to $X$
if $\lim_{n\to\infty}\|X_n-X\|_{\infty}=0$, which is denoted by
$X_n\stackrel{L^{\infty}}{\longrightarrow} X$.

\item $\{X_n,n\geq 1\}$  is said to converge to $X$ in distribution, if for any bounded continuous function $f$, $\lim_{n\to\infty}E[f(X_n)]=E[f(X)]$,  which is denoted by $X_n\stackrel{d}{\longrightarrow} X$.

\item $\{X_n,n\geq 1\}$  is said to completely converge to $X$, if for any $\varepsilon>0$, $\sum_{n=1}^{\infty}P(\{|X_n-X|\geq \varepsilon\})<\infty$,  which is denoted by $X_n\stackrel{c.c.}{\longrightarrow} X$ (see \cite{HR47}).

\item $\{X_n,n\geq 1\}$  is said to S-$L^p$ converge to $X$
$(p>0)$ if $\sum_{n=1}^{\infty}E[|X_n-X|^p]<\infty$, which is
denoted by $X_n\stackrel{S\mbox{-}L^p}{\longrightarrow} X$ (see
\cite[Definition 1.4]{LH15}).
\end{itemize}
The relations among the different kinds of convergence can be
described as follows.

\begin{eqnarray*}
\begin{array}{ccccccccc}
X_n\stackrel{S\mbox{-}L^p}{\longrightarrow} X &\Rightarrow & X_n\stackrel{c.c.}{\longrightarrow} X& \Rightarrow & X_n\stackrel{a.s.}{\longrightarrow} X &\Rightarrow & X_n\stackrel{P}{\longrightarrow} X &\Rightarrow & X_n\stackrel{d}{\longrightarrow} X,\\
&&&&\Uparrow&&\Uparrow&&\\
&&&&X_n\stackrel{L^{\infty}}{\longrightarrow} X &\Rightarrow&X_n\stackrel{L^p}{\longrightarrow} X&&
\end{array}
\end{eqnarray*}
and
\begin{itemize}

\item if $X_n\stackrel{P}{\longrightarrow} X$, then there exists a subsequence $\{X_{n_k}\}$ of $\{X_n\}$ such that $X_{n_k}\stackrel{a.s.}{\longrightarrow} X$ as $k\to\infty$;

\item if $X_n\stackrel{d}{\longrightarrow} C$, where $C$ is a constant, then $X_n\stackrel{P}{\longrightarrow} C$;

\item if $X_n\stackrel{d}{\longrightarrow} X$, then by  Skorokhod's theorem, there exists a sequence of random variables $\{Y,Y_n,n\geq 1\}$ such that for any $n\geq 1$, $X_n$ and $Y_n$ have the same distribution, $X$ and $Y$ have the same distribution, and $Y_n\stackrel{a.s.}{\longrightarrow} Y$.
\end{itemize}

In virtue of the relation between convergence in probability and
complete convergence, the relation between $L^p$ convergence and
S-$L^p$ convergence,  we  introduce two new kinds of convergence
of random variables, which are stronger versions of a.s.
convergence and $L^{\infty}$ convergence, respectively.

\begin{defi}\label{defi-1.1} Let $\alpha>0$.
$\{X_n,n\geq 1\}$  is said to strongly almost surely converge to
$X$ with order $\alpha$, if
$$
\sum_{n=1}^{\infty}|X_n-X|^{\alpha}<\infty\ a.s.,
$$
 which is denoted by $X_n\stackrel{S_{\alpha}\mbox{-}a.s.}{\longrightarrow} X$.
\end{defi}

\begin{defi}\label{defi-1.2}
$\{X_n,n\geq 1\}$  is said to strongly $L^{\infty}$-converge to
$X$ if
$$
\sum_{n=1}^{\infty}\|X_n-X\|_{\infty}<\infty,
$$
 which is denoted by $X_n\stackrel{S\mbox{-}L^{\infty}}{\longrightarrow} X$.
\end{defi}

We now introduce two new kinds of convergence which are stronger
versions of convergence in distribution.

\begin{defi}\label{defi-1.5}
$\{X_n,n\geq 1\}$  is said to $S_1\mbox{-}d$ converge to $X$, if
for any bounded Lipschitz continuous function $f$,
$$
\sum_{n=1}^{\infty}|E[f(X_n)-f(X)]|<\infty,
$$
which is denoted by $X_n\stackrel{S_1\mbox{-}d}{\longrightarrow}
X$.
\end{defi}

\begin{defi}\label{defi-1.8}
Let $F_n$ and $F$ be the distribution functions of $X_n$ and $X$,
respectively. $\{X_n,n\geq 1\}$  is said to $S_2\mbox{-}d$
converge to $X$, if for any continuous point $x$ of  $F$,
$$
\sum_{n=1}^{\infty}|F_n(x)-F(x)|<\infty,
$$
which is denoted by $X_n\stackrel{S_2\mbox{-}d}{\longrightarrow}
X$.
\end{defi}

The rest of this paper is organized as follows. In
Section 2, we study the law of large numbers for
independent and identically distributed (i.i.d.)  random
variables. In particular, we obtain a strong $L^p$-convergence
version and a strongly almost sure convergence version of the law
of large numbers. In Section 3, we discuss the relations among several kinds of convergence.
In Section 4, we present some open questions for further research.

\section{Convergence rates in the law of large numbers}\setcounter{equation}{0}

Let $\{X,X_n,n\geq 1\}$ be a sequence of i.i.d. random variables.
Define $S_n=X_1+\cdots+X_n$, $n\in\mathbb{N}$. Hsu and Robbins
(\cite{HR47}) proved that if $E[X^2]<\infty$ and $E[X]=\mu$, then
$\frac{S_n}{n}\stackrel{c.c.}{\longrightarrow} \mu$. Erd\"{o}s
(\cite{Er49}) proved the converse result. Baum and Katz
(\cite{BK65}) extended the Hsu-Robbins-Erd\"{o}s theorem. Below is
a special case of the Baum-Katz theorem.

\begin{thm} (Baum and Katz \cite{BK65}). Let $\alpha\geq 1$. Suppose that $\{X, X_n,n\geq 1\}$ is a sequence of i.i.d.
random variables with partial sum $S_n=\sum_{i=1}^n X_i$,
$n\in\mathbb{N}$. Then, the condition $E|X|^{\alpha}<\infty$ and
$EX=0$ is equivalent to
$$
\sum_{n=1}^{\infty}n^{\alpha-2}P(|S_n|>n\epsilon)<\infty,\ \
\forall\epsilon>0.
$$
\end{thm}
Lanzinger (\cite{La98}), Gut and Stadtm\"{u}ller (\cite{GS11}),
Chen and Sung (\cite{CS14}) extended the results of Baum and Katz.

Chow (\cite{Ch88}) first investigated the complete moment
convergence and obtained the following result. Let $\alpha\geq 1$,
$p\le \alpha$ and $p<2$. Suppose that $\{X, X_n,n\geq 1\}$ is a
sequence of i.i.d. random variables with $E[X]=0$. If
$E[|X|^{\alpha}+|X|\log^+|X|]<\infty$, then
\begin{eqnarray*}\label{Chow-a}
\sum_{n=1}^{\infty}n^{\frac{\alpha}{p}-\frac{1}{p}-2}E\left[\left(|S_n|-\varepsilon
n^{\frac{1}{p}}\right)^{+}\right]<\infty\ \mbox{for all}\
\varepsilon>0,
\end{eqnarray*}
where $x^+=\max\{0,x\}$.

Chow's result has been generalized in various directions. Wang and Su (\cite{WS04}), Wang et al. (\cite{WZW05}), Chen (\cite{Ch06}),  Guo and Xu (\cite{GX06}), Rosalsky et al. (\cite{RTV06}), Ye and Zhu (\cite{YZ07}), and Qiu et al. (\cite{QC14}) studied the
complete moment convergence for sums of Banach space valued random elements. Li and Zhang (\cite{LZ04}), Chen et al. (\cite{CHV07}), Kim et al. (\cite{KKC08}), and Zhou (\cite{Zh10}) considered the complete moment convergence for moving average processes. Jiang and Zhang (\cite{JZ06}), Li (\cite{Li06}), Liu and Lin (\cite{LL06}), Ye et al. (\cite{YZP07}), Fu and Zhang (\cite{FZ08}), Zhao and Tao (\cite{ZT08}), and Chen and Zhang (\cite{CZ10}) studied precise asymptotics for the complete moment convergence. Wang and Zhao (\cite{WZ06}), Liang et al. (\cite{LLR10}), and Guo (\cite{Gu13}) considered the complete moment convergence for negatively associated random variables. Qiu and Chen (\cite{QC14}) studied the complete moment convergence for i.i.d. random variables and extended two results in Gut and Stadtm\"{u}ller (\cite{GS11}) to the complete moment convergence.

In the following of this section, we  study convergence rates in the law of large numbers for
i.i.d. random variables. In
particular, we obtain a strong $L^p$-convergence version and a
strongly almost sure convergence version of the law of large
numbers.

\subsection{Strong $L^p $-convergence version of the law of large numbers}

Let $\{Y,Y_n,n\geq 1\}$ be a sequence of random variables and
$p>0$. We have
\begin{eqnarray}\label{3.1}
Y_n\stackrel{S\mbox{-}L^p}{\longrightarrow} Y \Rightarrow
Y_n\stackrel{c.c.}{\longrightarrow} Y.
\end{eqnarray}
In this subsection, we consider the following question:

{\it Does it hold that
$\frac{S_n}{n}\stackrel{S\mbox{-}L^p}{\longrightarrow} \mu$ for
some $p>0$ under some condition? }

By the Hsu-Robbins-Erd\"{o}s theorem and (\ref{3.1}),  we know
that the condition $E[X^2]<\infty$ is needed in order that
$\frac{S_n}{n}\stackrel{S\mbox{-}L^p}{\longrightarrow} \mu$.

\begin{thm}\label{new1}
(1) If $E[X^2]<\infty$ and $X\not\equiv\mu$ a.s., then
$\frac{S_n}{n}\stackrel{S\mbox{-}L^p}{\nrightarrow} \mu$ for any
$0<p\le 2$.

\noindent (2) If $\alpha>2$ and $E[|X|^{\alpha}]<\infty$, then
$\frac{S_n}{n}\stackrel{S\mbox{-}L^p}{\rightarrow} \mu$ for any
$2<p\le \alpha$.
\end{thm}
{\bf Proof.} We assume without loss of generality that $\mu=0$ and
$E[X^2]=1$.

\noindent (1) We have
$$
E\left[\left|\frac{S_n}{n}\right|^{p}\right]=
\frac{1}{n^{p/2}}E\left[\left|\frac{S_n}{\sqrt{n}}\right|^{p}\right].
$$
Denote the distribution function of $\frac{S_n}{\sqrt{n}}$ by
$F_n$. Let $f\in C_c(\mathbb{R})$ satisfying $|f(x)|= |x|^{p}$ for
$|x|\le 1$ and $|f(x)|\le |x|^{p}$ for $|x|> 1$. Then, by the
central limit theorem, we have
$$
E\left[\left|\frac{S_n}{\sqrt{n}}\right|^{p}\right]\ge\int_{-\infty}^{\infty}|f(x)|dF_n(x)\rightarrow
\int_{-\infty}^{\infty}|f(x)|\frac{1}{\sqrt{2\pi}}e^{-|x|^2/2}dx>0.
$$
Define
$$
c=\int_{-\infty}^{\infty}|f(x)|\frac{1}{\sqrt{2\pi}}e^{-|x|^2/2}dx.
$$
Then, there exists $N\in \mathbb{N}$ such that
$$
E\left[\left|\frac{S_n}{\sqrt{n}}\right|^{p}\right]\ge
\frac{c}{2},\ \ \forall n\ge N.
$$
Therefore,
$$
\sum_{n=1}^{\infty}E\left[\left|\frac{S_n}{n}\right|^{p}\right]
\ge\sum_{n=N}^{\infty}\frac{1}{n^{p/2}}E\left[\left|\frac{S_n}{\sqrt{n}}\right|^{p}\right]\ge
\frac{c}{2}\sum_{n=N}^{\infty}\frac{1}{n^{p/2}}=\infty.
$$

\noindent (2) If $Y$ is a random variable, we denote
$\|Y\|_{L^r}:=(E[|Y|^r])^{1/r}$ for $r\ge 1$. By the
Burkholder-Davis-Gundy inequality and Minkowski's inequality, we
have
\begin{eqnarray}\label{p1}
E[|S_n|^{\alpha}]&\le& cE[(X_1^2+\cdots +X_n^2)^{\alpha/2}]\nonumber\\
&=&c\|X_1^2+\cdots +X_n^2\|^{\alpha/2}_{L^{\alpha/2}}\nonumber\\
&\le&c(\|X_1^2\|_{L^{\alpha/2}}+\cdots +\|X_n^2\|_{L^{\alpha/2}})^{\alpha/2}\nonumber\\
&=&cn^{\alpha/2}E[|X|^{\alpha}],
\end{eqnarray}
where $c>0$ is a constant, which is independent of $n$. Then,
$$
\sum_{n=1}^{\infty}E\left[\left|\frac{S_n}{n}\right|^{\alpha}\right]
\le c\sum_{n=1}^{\infty}n^{-\alpha/2}E[|X|^{\alpha}]<\infty.
$$

For $2<p< \alpha$, we obtain by (\ref{p1}) that
\begin{eqnarray*}
\sum_{n=1}^{\infty}E\left[\left|\frac{S_n}{n}\right|^{p}\right]&\le& \sum_{n=1}^{\infty}\frac{1}{n^p}\left(E\left[\left|S_n\right|^{\alpha}\right]\right)^{p/\alpha}\\
&\le&\sum_{n=1}^{\infty}\frac{1}{n^p}(cn^{\alpha/2}E[|X|^{\alpha}])^{p/\alpha}\\
&=&\sum_{n=1}^{\infty}\frac{1}{n^{p/2}}(c
E[|X|^{\alpha}])^{p/\alpha}\\
&<&\infty.
\end{eqnarray*}
\hfill\fbox

In \cite{Ch88},  Chow also obtained the following result. Let $\{X,X_n,n\geq 1\}$ be a sequence of i.i.d. random variables with $E[X]=0$. Suppose that $1<\alpha<2$. If $E[|X|^{\alpha}\log^{+}|X|]<\infty$, then
\begin{eqnarray*}\label{Chow-a}
\sum_{n=1}^{\infty}n^{-2}E[|S_n|^{\alpha}]<\infty.
\end{eqnarray*}
As a direct consequence of Theorem \ref{new1} and its proof, we
have the following corollaries.

\begin{cor}
Suppose that $\alpha>2$, $E[|X|^{\alpha}]<\infty$ and $E[X]=0$.
Then, for any $2<p\le\alpha$ and $\beta>(p+2)/2$, we have
$$
\sum_{n=1}^{\infty}n^{-\beta}E[|S_n|^p]<\infty.
$$
\end{cor}

\begin{cor}
Suppose that $X\not\equiv\mu$ a.s. and  $E[|X|^{\alpha}]<\infty$
for any $\alpha>0$. Then
$\frac{S_n}{n}\stackrel{S\mbox{-}L^{p}}{\longrightarrow} \mu$ if
and only if $p>2$.
\end{cor}

%\begin{exa}
%Suppose that $X\sim N(\mu,\sigma^2) (\sigma>0)$. Then $\frac{S_n}{n}-\mu\sim N(0,\frac{\sigma^2}{n})$ and so
%$\frac{S_n-n\mu}{\sqrt{n}}\sim N(0,\sigma^2)$ .  For any $\alpha>0$, we have
%\begin{eqnarray*}
%E\left[\left|\frac{S_n}{n}-\mu\right|^{\alpha}\right]&=&
%\frac{1}{\sqrt{n}^{\alpha}}E\left[\left|\frac{S_n-n\mu}{\sqrt{n}}\right|^{\alpha}\right]=c\frac{1}{\sqrt{n}^{\alpha}},
%\end{eqnarray*}
%where $c=\int_{-\infty}^{\infty}|x|^{\alpha}\frac{1}{\sqrt{2\pi}\sigma}e^{-\frac{x^2}{2\sigma^2}}dx>0$.  It follows that the series
%$
%\sum_{n=1}^{\infty}E\left[\left|\frac{S_n}{n}-\mu\right|^{\alpha}\right]
%$
%converges if and only if $\alpha>2$, i.e. in this case $\frac{S_n}{n}\stackrel{S\mbox{-}L^{\alpha}}{\longrightarrow} \mu$ if and only if $\alpha>2$.
%\end{exa}

%\begin{rem}
%By Propositions \ref{pro-2.1} and \ref{new1}, we know that if
%$X\not\equiv\mu$ a.s., then
%$\frac{S_n}{n}\stackrel{S\mbox{-}L^{\infty}}{\nrightarrow} \mu$.
%In fact, it is easy to know that if $X\not\equiv\mu$ a.s., then
%$\frac{S_n}{n}\stackrel{L^{\infty}}{\nrightarrow} \mu$.
%\end{rem}

\subsection{Strongly almost sure convergence version of the law of large numbers}

In this subsection, we consider the following question:

{\it Does it hold that
$\frac{S_n}{n}\stackrel{S_{\alpha}\mbox{-}a.s.}{\longrightarrow}
\mu$ for some $\alpha>0$ under some condition? }

\begin{thm}\label{april3}
(1) If $E[X^4]<\infty$ and $X\not\equiv\mu$ a.s., then
$\sum_{n=1}^{\infty}\left|\frac{S_n}{n}-\mu\right|^{\alpha}=\infty$ a.s.
for any $0<\alpha\le 2$.

\noindent (2) If $E[|X|^2]<\infty$, then for any $\alpha>2$ we
have
$$
\frac{S_n}{n}\stackrel{S_{\alpha}\mbox{-}a.s.}{\longrightarrow}
\mu.
$$
\end{thm}
{\bf Proof.} We assume without loss of generality that $\mu=0$ and
$E[X^2]=1$.

\noindent (1) For $N\in\mathbb{N}$, we have
\begin{eqnarray*}
E\left[\left|\frac{S_n}{n}\right|^2\right]=\frac{1}{n^2}E\left[\sum_{i=1}^n X_i^2+2\sum_{1\leq i<j\leq n}X_iX_j\right]=\frac{1}{n}.
\end{eqnarray*}
Define
\begin{eqnarray*}
e_N=\sum_{n=1}^{N}\frac{1}{n},\quad\quad
W_N=e_N+2\sum_{n=2}^{N}\frac{\sum_{1\le i<j\le n}X_iX_j}{n^2},
\end{eqnarray*}
and
$$
R_N=\sum_{n=1}^{N}\left(\frac{X_1^2+\cdots+ X_n^2}{n^2}-\frac{1}{n}\right).
$$ Then,
\begin{equation}\label{swim}
\sum_{n=1}^{N}\left|\frac{S_n}{n}\right|^2=W_N+R_N.
\end{equation}

For $M,N\in\mathbb{N}$ with $M<N$, we have
\begin{eqnarray*}
&&E[(R_N-R_M)^2]\\
&=&E\left[\left(\sum_{n=M+1}^{N}\frac{(X_1^2-1)+\cdots+ (X_n^2-1)}{n^2}\right)^2\right]\\
&=&\sum_{n=M+1}^{N}\frac{E[(X^2-1)^2]}{n^3}+2\sum_{M+1\le k<l\le N}\frac{E[(\sum_{i=1}^k(X_i^2-1))(\sum_{j=1}^l(X_j^2-1))]}{k^2l^2}\\
&=&\sum_{n=M+1}^{N}\frac{E[(X^2-1)^2]}{n^3}+2\sum_{M+1\le k<l\le N}\frac{E[(X^2-1)^2]}{kl^2}\\
&\le&(E[(X^2-1)^2])\left\{\sum_{n=M+1}^{N}\frac{1}{n^3}+2\sum_{l=M+1}^N\frac{e_l}{l^2}\right\}\\
&\rightarrow&0\ \ \ \ \ \ {\rm as}\ \ \ M\rightarrow\infty.
\end{eqnarray*}
Hence $\{R_N\}_{N=1}^{\infty}$ is a Cauchy sequence in $L^2$,
which implies that $\{R_N\}_{N=1}^{\infty}$ converges to some
$R\in L^2$ in probability.

Suppose that $\frac{S_n}{n}\stackrel{S_2\mbox{-}a.s.}{\rightarrow}
0$. Dividing both sides of  (\ref{swim}) by $e_N$ and letting
$N\rightarrow\infty$, we get
$$
1+\frac{2}{e_N}\sum_{n=2}^{N}\frac{\sum_{1\le i<j\le n}X_iX_j}{n^2}\rightarrow 0\ \ {\rm in\ probability}\ \ {\rm as}\ N\rightarrow\infty,
$$
which implies that
\begin{equation}\label{con1}
\frac{1}{e_N}\sum_{n=2}^{N}\frac{\sum_{1\le i<j\le n}X_iX_j}{n^2}\rightarrow -\frac{1}{2}\ \ {\rm in\ probability}\ \ {\rm as}\ N\rightarrow\infty.
\end{equation}
We have
\begin{eqnarray}\label{con3}
&&E\left[\left(\frac{1}{e_N}\sum_{n=2}^{N}\frac{\sum_{1\le i<j\le n}X_iX_j}{n^2}\right)^2\right]\nonumber\\
&=&\frac{1}{(e_N)^2}\left(\sum_{n=2}^N\frac{E\left[(\sum_{1\leq i<j\leq n}X_iX_j)^2\right]}{n^4}\right.\nonumber\\
&&\left.+2\sum_{2\leq k<l\leq N}\frac{1}{k^2l^2}E\left[\left(\sum_{1\leq i<j\leq k}X_iX_j\right)\left(\sum_{1\leq i<j\leq l}X_iX_j\right)\right]\right)\nonumber\\
&=&\frac{1}{(e_N)^2}\left(\sum_{n=2}^{N}\frac{n-1}{2n^3}+\sum_{2\le k<l\le N}\frac{k-1}{kl^2}\right)\nonumber\\
&\le&\frac{1}{(e_N)^2}\left(\sum_{n=2}^{N}\frac{1}{n^2}+e_N\right)\nonumber\\
&\le&\sum_{n=1}^{\infty}\frac{1}{n^2}.
\end{eqnarray}
By (\ref{con1}) and (\ref{con3}), we get
$$
\lim_{N\rightarrow\infty}E\left[\frac{1}{e_N}\sum_{n=2}^{N}\frac{\sum_{1\le i<j\le n}X_iX_j}{n^2}\right]=-\frac{1}{2},
$$
which contradicts with
$$
E\left[\frac{1}{e_N}\sum_{n=2}^{N}\frac{\sum_{1\le i<j\le n}X_iX_j}{n^2}\right]=0, \ \ \forall N\ge2.
$$
Then, $\frac{S_n}{n}\stackrel{S_2\mbox{-}a.s.}{\nrightarrow}
0$.  Therfore, we obtain by the  Hewitt-Savage  0-1 law that $\sum_{n=1}^{\infty}\left|\frac{S_n}{n}\right|^2=\infty$ a.s..

By the strong law of large numbers, there exists a set $N\in
\mathcal{F}$ satisfying $P(N)=0$ and for any $\omega\in
\Omega\backslash N$, there exists $M(\omega)\in\mathbb{N}$ such
that for any $n\ge M(\omega)$,
$$
\frac{|S_n(\omega)|}{n}<1. $$
 It follows that for any $0<\alpha<2$ and
$\omega\in \Omega\backslash N$,
$$
\sum_{n=M(\omega)}^{\infty}\left|\frac{S_n(\omega)}{n}\right|^2\le
\sum_{n=M(\omega)}^{\infty}\left|\frac{S_n(\omega)}{n}\right|^{\alpha}.
$$
Therefore, $\sum_{n=1}^{\infty}\left|\frac{S_n}{n}\right|^2=\infty$ a.s. implies that
$\sum_{n=1}^{\infty}\left|\frac{S_n}{n}\right|^{\alpha}=\infty$ a.s.
for any $0<\alpha< 2$.

\noindent (2) By the Hartman-Wintner law of iterated logarithm, we
have
$$
\limsup_{n\to\infty}\frac{|S_n|}{\sqrt{2n \log\log n}}=1\ a.s..
$$
Then, there exists a set $N\in \mathcal{F}$ satisfying $P(N)=0$
and for any $\omega\in \Omega\backslash N$, there exists
$M(\omega)\in\mathbb{N}$ such that for any $n\ge M(\omega)$,
$$
\frac{|S_n(\omega)|}{\sqrt{2n \log\log n}}<2.
$$
It follows that for any $\alpha>2$ and $\omega\in \Omega\backslash N$,
$$
\sum_{n=1}^{\infty}\left|\frac{S_n}{n}\right|^{\alpha}<\infty,
$$
i.e., $\frac{S_n}{n}\stackrel{S_{\alpha}\mbox{-}a.s.}{\longrightarrow}0$.\hfill\fbox

\begin{rem}
By analogues of the Hartman-Wintner law of iterated logarithm in
the infinite variance case, we can show that
$\frac{S_n}{n}\stackrel{S_{\alpha}\mbox{-}a.s.}{\longrightarrow}
\mu$ for any $\alpha>2$ under weaker conditions. Define $L(x) =
\log\max\{e,x\}$ and $LL(x) = L(Lx)$ for $x\in\mathbb{R}$. By
Einmahl and Li \cite[Corollaries 1 and 2]{EL}, we have
$\frac{S_n}{n}\stackrel{S_{\alpha}\mbox{-}a.s.}{\longrightarrow}
\mu$ for any $\alpha>2$ if one of the following conditions is
fulfilled.

(i) For some $p\ge1$,
$$
E\left[\frac{(X-\mu)^2}{(LL(|X-\mu|))^p}\right]<\infty,\ \
\limsup_{x\rightarrow\infty}(LL(x))^{1-p}E[(X-\mu)^21_{\{|X-\mu|\le
x\}}]<\infty.
$$

(ii) For some $r>0$,
$$
E\left[\frac{(X-\mu)^2}{(L(|X-\mu|))^r}\right]<\infty,\ \
\limsup_{x\rightarrow\infty}\frac{LL(x)}{(L(x))^r}E[(X-\mu)^21_{\{|X-\mu|\le
x\}}]<\infty.
$$
\end{rem}

\section{Relations among several kinds of convergence}\setcounter{equation}{0}

\subsection{Main results}

\begin{pro}\label{pro-2.1}

Let $\{X,X_n,n\geq 1\}$ be a sequence of random variables. Then

(i)  For $p\geq 1$, we have
$X_n\stackrel{S\mbox{-}L^{\infty}}{\longrightarrow} X\Rightarrow
X_n\stackrel{S\mbox{-}L^p}{\longrightarrow} X$;

%(ii)  For $r\leq 1$, we have $X_n\stackrel{S\mbox{-}L^r}{\longrightarrow} X\Rightarrow X_n\stackrel{S\mbox{-}L^1}{\longrightarrow} X$.

(ii) For $\alpha>0$, we have
$X_n\stackrel{S\mbox{-}L^{\alpha}}{\longrightarrow} X\Rightarrow
X_n\stackrel{S_{\alpha}\mbox{-}a.s.}{\longrightarrow} X$;

(iii) For any $\alpha\ge1$, we have
$X_n\stackrel{S\mbox{-}L^{\infty}}{\longrightarrow} X\Rightarrow
X_n\stackrel{S_{\alpha}\mbox{-}a.s.}{\longrightarrow} X$.
\end{pro}

It is well known that $X_n\stackrel{P}{\longrightarrow} X$ if and
only if for any subsequence $\{X_{n'}\}$ of $\{X_n\}$, there
exists a subsequence $\{X_{n'_k}\}$ of  $\{X_{n'}\}$  such that
$X_{n'_k}\stackrel{a.s.}{\longrightarrow} X$.  In the following,
we strengthen this result to the strongly almost sure convergence.

\begin{thm}\label{thm-2.2}
 $X_n\stackrel{P}{\longrightarrow} X$ if and only if for any subsequence $\{X_{n'}\}$ of $\{X_n\}$ and some (hence all) $\alpha>0$, there exists a subsequence
 $\{X_{n'_k}\}$ of  $\{X_{n'}\}$  such that $X_{n'_k}\stackrel{S_{\alpha}\mbox{-}a.s.}{\longrightarrow} X$.
\end{thm}

As a direct consequence of Theorem \ref{thm-2.2}, we have the
following corollary.

\begin{cor}\label{cor-2.3}
If  $X_n\stackrel{c.c.}{\longrightarrow} X$, or
$X_n\stackrel{a.s.}{\longrightarrow} X$, or
$X_n\stackrel{L^{\infty}}{\longrightarrow} X$, or
$X_n\stackrel{L^p}{\longrightarrow} X$, then for any $\alpha>0$,
there is a subsequence $\{X_{n_k}\}$ of $\{X_n\}$  such that
$X_{n_k}\stackrel{S_{\alpha}\mbox{-}a.s.}{\longrightarrow} X$.
\end{cor}

\begin{pro}\label{pro-2.4}
$X_n\stackrel{S\mbox{-}L^1}{\longrightarrow} X\Rightarrow
X_n\stackrel{S_1\mbox{-}d}{\longrightarrow} X.$
\end{pro}

\bigskip

Now we have the following diagram:
\begin{eqnarray*}
\begin{array}{ccccccc}
&& X_n\stackrel{S_1\mbox{-}d}{\longrightarrow} X & \Rightarrow  & \Rightarrow & \Rightarrow & X_n\stackrel{d}{\longrightarrow} X\\
&&\Uparrow & &&&\Uparrow\\
X_n\stackrel{S\mbox{-}L^{\infty}}{\longrightarrow} X & \Rightarrow & X_n\stackrel{S\mbox{-}L^1}{\longrightarrow} X & \Rightarrow &
X_n\stackrel{S_1\mbox{-}a.s.}{\longrightarrow} X&&\Uparrow\\
&&\Downarrow&&\Downarrow &&\Uparrow\\
&&X_n\stackrel{c.c.}{\longrightarrow} X& \Rightarrow& X_n\stackrel{a.s.}{\longrightarrow} X& \Rightarrow& X_n\stackrel{P}{\longrightarrow} X.\\
&&&&\Uparrow& & \Uparrow \\
&&&&X_n\stackrel{L^{\infty}}{\longrightarrow} X&\Rightarrow & X_n\stackrel{L^1}{\longrightarrow} X
\end{array}
\end{eqnarray*}

\begin{thm}\label{pro-2.5}
Let $C$ be a constant. Then
$X_n\stackrel{S_2\mbox{-}d}{\longrightarrow} C\Leftrightarrow
X_n\stackrel{c.c.}{\longrightarrow} C$.
\end{thm}

If $C$ is a constant and
$Y_n\stackrel{S_2\mbox{-}d}{\longrightarrow} C$, then by Theorem
\ref{pro-2.5} we know that $Y_n\stackrel{c.c.}{\longrightarrow}
C$. In the following Slutsky-type theorem, we need a  stronger
condition than $Y_n\stackrel{c.c.}{\longrightarrow} C$.

\begin{pro}\label{pro-2.6}
Suppose that $X_n\stackrel{S_1\mbox{-}d}{\longrightarrow} X$ and
$Y_n\stackrel{S\mbox{-}L^1}{\longrightarrow} C$. Then

(i) $X_n+Y_n\stackrel{S_1\mbox{-}d}{\longrightarrow} X+C$;

(ii) if $\{X_n\}$ is a sequence of bounded random variables, then
$X_nY_n\stackrel{S_1\mbox{-}d}{\longrightarrow} CX$;

(iii) if  $\{X_n\}$ and $\{1/Y_n\}$  are two  sequences of bounded
random variables and $C\neq 0$, then
$\frac{X_n}{Y_n}\stackrel{S_1\mbox{-}d}{\longrightarrow}
\frac{X}{C}$.
\end{pro}

\begin{pro}\label{pro-2.7}
Let $\{X,X_n,n\geq 1\}$ be a sequence of random variables and
$\{F,F_n,n\geq 1\}$ be the corresponding sequence of distribution
functions. Then $X_n\stackrel{S_2\mbox{-}d}{\longrightarrow} X$ if
one of the following conditions is fulfilled.

(1) $X$ is a discrete random variable such that $\{x\in\mathbb{R}:P(X=x)=0\}$ is an open subset of
$\mathbb{R}$ and $X_n\stackrel{c.c.}{\longrightarrow} X$.

(2) $X$ has a bounded density function and
$\sum_{n=1}^{\infty}P\{n(\log n)^{1+\beta}|X_n-X|\geq
\delta\}<\infty$ for two positive constants $\beta$ and $\delta$.
\end{pro}

\subsection{Proofs}

 {\bf Proof of Proposition \ref{pro-2.1}.}  
 
 (i) If $\|X_n-X\|_{\infty}<1$ and $p\ge 1$, we have
$$
E[|X_n-X|^p]\leq \|X_n-X\|_{\infty}^p\leq \|X_n-X\|_{\infty},
$$
which together with the definitions of
$X_n\stackrel{S\mbox{-}L^{\infty}}{\longrightarrow} X$ and
$X_n\stackrel{S\mbox{-}L^p}{\longrightarrow} X$ implies that
$X_n\stackrel{S\mbox{-}L^{\infty}}{\longrightarrow} X\Rightarrow
X_n\stackrel{S\mbox{-}L^p}{\longrightarrow} X$.

(ii) Let $\alpha>0$. If
$X_n\stackrel{S\mbox{-}L^{\alpha}}{\longrightarrow} X$, then
$\sum_{n=1}^{\infty}E[|X_n-X|^{\alpha}]<\infty$. By the monotone
convergence theorem, we have
$$
\sum_{n=1}^{\infty}E[|X_n-X|^{\alpha}]=E\left[\sum_{n=1}^{\infty}|X_n-X|^{\alpha}\right].
$$
It follows that $E[\sum_{n=1}^{\infty}|X_n-X|^{\alpha}]<\infty$
and thus $\sum_{n=1}^{\infty}|X_n-X|^{\alpha}<\infty$ a.s., i.e.,
$X_n\stackrel{S_{\alpha}\mbox{-}a.s.}{\longrightarrow} X$.

(iii) It is a direct consequence of (i) and (ii). \hfill\fbox\\

\noindent {\bf Proof of Theorem \ref{thm-2.2}.} The sufficiency is
obvious. We only prove the necessity. Suppose that
$X_n\stackrel{P}{\longrightarrow} X$ and $\{X_{n'}\}$ is a
subsequence of $\{X_n\}$. Then
$X_{n'}\stackrel{P}{\longrightarrow} X$. Thus, for any $k\in
\mathbb{N}$, we have
$$
\lim_{n'\to\infty}P\left\{|X_{n'}-X|^{\alpha}\geq
\frac{1}{k^2}\right\}=0.
$$
It follows that there exists a sequence $\{X_{n'_k}\}$ of  $\{X_{n'}\}$ such that for any $k\in \mathbb{N}$,
$$
P\left\{|X_{n'_k}-X|^{\alpha}\geq \frac{1}{k^2}\right\}\leq
\frac{1}{k^2},
$$
which implies that
$$
\sum_{k=1}^{\infty}P\left\{|X_{n'_k}-X|^{\alpha}\geq
\frac{1}{k^2}\right\}<\infty.
$$
By the Borel-Cantelli lemma, we get
$$
P\left(\bigcap_{n=1}^{\infty}\bigcup_{k=n}^{\infty}\left\{|X_{n'_k}-X|^{\alpha}\geq
\frac{1}{k^2}\right\}\right)=0,
$$
which implies that
$$
P\left(\bigcup_{n=1}^{\infty}\bigcap_{k=n}^{\infty}\left\{|X_{n'_k}-X|^{\alpha}<
\frac{1}{k^2}\right\}\right)=1.
$$
Therefore, $X_{n'_k}\stackrel{S_{\alpha}\mbox{-}a.s.}{\longrightarrow} X$.  The proof is complete. \hfill\fbox\\

\noindent {\bf Proof of Proposition \ref{pro-2.4}.} Suppose that
$f$ is a bounded Lipschitz continuous function. Then there a
positive constant $C$ such that
$$
|f(x)-f(y)|\leq C|x-y|,\ \ \forall x,y\in\mathbb{R}.
$$
It follows that
$$
|E[f(X_n)-f(X)]|\leq E[|f(X_n)-f(X)|]\leq CE[|X_n-X|],
$$
which together with the definitions of
$X_n\stackrel{S\mbox{-}L^1}{\longrightarrow} X$ and
$X_n\stackrel{S_1\mbox{-}d}{\longrightarrow} X$ implies that
$X_n\stackrel{S\mbox{-}L^1}{\longrightarrow} X\Rightarrow X_n\stackrel{S_1\mbox{-}d}{\longrightarrow} X$.\hfill\fbox\\

\noindent {\bf Proof of Theorem \ref{pro-2.5}.} 

``$\Rightarrow$"
For any $\epsilon>0$, we have
\begin{eqnarray}\label{s_6-cc-1}
P\{|X_n-C|\geq \epsilon\}&=&1-P\{X_n<C+\epsilon\}+P\{X_n\le
C-\epsilon\}\nonumber\\
&\leq& 1-F_n\left(C+\frac{\epsilon}{2}\right)+F_n(C-\epsilon).
\end{eqnarray}
If $X_n\stackrel{S_2\mbox{-}d}{\longrightarrow} C$, then for any
$\epsilon>0$,
\begin{eqnarray*}
&&\sum_{n=1}^{\infty}\left|F_n\left(C+\frac{\epsilon}{2}\right)-1\right|<\infty\
\ \ \mbox{and}\ \
\sum_{n=1}^{\infty}|F_n\left(C-\epsilon\right)-0|<\infty,
\end{eqnarray*}
which together with (\ref{s_6-cc-1}) implies that for any
$\epsilon>0$,
$$
\sum_{n=1}^{\infty}P\{|X_n-C|\geq \epsilon\}<\infty.
$$
Hence $X_n\stackrel{c.c.}{\longrightarrow} C$.

``$\Leftarrow$" Let $X\equiv C$ and denote by $F$ the distribution
of $X$. For any $\epsilon>0$ and $x\in \mathbb{R}$, we have
\begin{eqnarray}\label{re}
F(x-\epsilon)-P\{|X_n-X|\geq \epsilon\}\leq F_n(x)\leq
P\{|X_n-X|\geq \epsilon\}+F(x+\epsilon).
\end{eqnarray}
If $x>C$, set $\epsilon=(x-C)/2$. By (\ref{re}), we have
\begin{eqnarray*}\label{s_6-cc-2}
1-P\{|X_n-C|\geq \epsilon\}\leq F_n(x)\leq P\{|X_n-C|\geq
\epsilon\}+1,
\end{eqnarray*}
i.e.,
\begin{eqnarray}\label{s_6-cc-2}
-P\{|X_n-C|\geq \epsilon\}\leq F_n(x)-1\leq P\{|X_n-C|\geq
\epsilon\}.
\end{eqnarray}
If $x<C$, set $\epsilon=(C-x)/2$. By (\ref{re}), we have
\begin{eqnarray*}
0-P\{|X_n-C|\geq \epsilon\}\leq F_n(x)\leq P\{|X_n-C|\geq
\epsilon\}+0,
\end{eqnarray*}
i.e.,
\begin{eqnarray}\label{s_6-cc-3}
-P\{|X_n-C|\geq \epsilon\}\leq F_n(x)-0\leq P\{|X_n-C|\geq
\epsilon\}.
\end{eqnarray}

By $X_n\stackrel{c.c.}{\longrightarrow} C$, (\ref{s_6-cc-2}) and
(\ref{s_6-cc-3}), we obtain that for any $x\neq C$,
$$
\sum_{n=1}^{\infty}|F_n(x)-F(x)|<\infty.
$$
Hence $X_n\stackrel{S_2\mbox{-}d}{\longrightarrow} C$.\hfill\fbox\\

\noindent {\bf Proof of Proposition \ref{pro-2.6}.}  Suppose that $f$ is a
bounded Lipschitz continuous function. Then there exits a positive
constant $K$ such that
$$
|f(x)-f(y)|\leq K|x-y|,\ \ \forall x,y\in\mathbb{R}.
$$

(i)  We have
\begin{eqnarray}\label{Slutsky-1}
&&|E[f(X_n+Y_n)]-E[f(X+C)]|\nonumber\\
&&\leq|E[f(X_n+Y_n)]-E[f(X_n+C)]|+|E[f(X_n+C)]-E[f(X+C)]|\nonumber\\
&&\leq KE[|Y_n-C|]+|E[f(X_n+C)]-E[f(X+C)]|.
\end{eqnarray}
By the assumption that
$Y_n\stackrel{S\mbox{-}L^1}{\longrightarrow} C$, we have
\begin{eqnarray}\label{Slutsky-2}
\sum_{n=1}^{\infty}E[|Y_n-C|]<\infty.
\end{eqnarray}
Define $g(x)=f(x+c)$. Then, $g$ is a bounded Lipschitz continuous
function and thus by the assumption that
$X_n\stackrel{S_1\mbox{-}d}{\longrightarrow} X$, we get
\begin{eqnarray}\label{Slutsky-3}
\sum_{n=1}^{\infty}|E[f(X_n+C)]-E[f(X+C)]|<\infty.
\end{eqnarray}
By (\ref{Slutsky-1})-(\ref{Slutsky-3}) and the definition of
$S_1\mbox{-}d$ convergence, we obtain that
$X_n+Y_n\stackrel{S_1\mbox{-}d}{\longrightarrow} X+C$.

 The proofs for (ii) and (iii)  are similar, so we omit the details.\hfill\fbox\\

\noindent {\bf Proof of Proposition \ref{pro-2.7}.}

(1) Suppose that $x\in \mathbb{R}$ with $P(X=x)=0$. Then, by the assumption that $\{x\in\mathbb{R}:P(X=x)=0\}$ is an open subset of
$\mathbb{R}$, there
exists $\epsilon>0$ such that
$$
F(x)=F(x+\epsilon)=F(x-\epsilon),
$$
which together with (\ref{re}) implies that
\begin{eqnarray*}
|F_n(x)-F(x)|&\leq& P\{|X_n-X|\geq \epsilon\}+|F(x+\epsilon)-F(x)|+|F(x-\epsilon)-F(x)|\nonumber\\
&=&P\{|X_n-X|\geq \epsilon\}.
\end{eqnarray*}
It follows that
$$
\sum_{n=1}^{\infty}|F_n(x)-F(x)|\le \sum_{n=1}^{\infty}P\{|X_n-X|\geq \epsilon\},
$$
which together with the definitions of complete convergence and
$S_2\mbox{-}d$ implies that
$X_n\stackrel{S_2\mbox{-}d}{\longrightarrow} X$.

(2) By the assumption, we know that there exists a positive
constant $C$ such that $|f(x)|\leq C,\forall x\in \mathbb{R}$. It
follows that for any $x,y\in \mathbb{R}$,
\begin{eqnarray}\label{cc-s-6-d-1}
|F(x)-F(y)|=\left|\int_x^yf(u)du\right|\leq C|y-x|.
\end{eqnarray}

By (\ref{re}) and (\ref{cc-s-6-d-1}), we have
\begin{eqnarray*}
|F_n(x)-F(x)|&\leq& P\{|X_n-X|\geq \epsilon\}+|F(x+\epsilon)-F(x)|+|F(x-\epsilon)-F(x)|\nonumber\\
&\leq&P\{|X_n-X|\geq \epsilon\}+2C\epsilon.
\end{eqnarray*}
Then,
\begin{eqnarray*}
\sum_{n=1}^{\infty}|F_n(x)-F(x)|&\leq& \sum_{n=1}^{\infty}\left(P\left\{|X_n-X|\geq \frac{\delta}{n(\log n)^{1+\beta}}\right\}+2C\frac{\delta}{n(\log n)^{1+\beta}}\right)\\
&=&\sum_{n=1}^{\infty}P\left\{n(\log n)^{1+\beta}|X_n-X|\geq
\delta\right\}+2C\delta\sum_{n=1}^{\infty}\frac{1}{n(\log
n)^{1+\beta}}\\
&<&\infty,
\end{eqnarray*}
and thus
$X_n\stackrel{S_2\mbox{-}d}{\longrightarrow} X$. \hfill\fbox

\subsection{Examples and remarks}

The following example shows that
$X_n\stackrel{S\mbox{-}L^{\infty}}{\longrightarrow} X$ is stronger
than $X_n\stackrel{S\mbox{-}L^p}{\longrightarrow} X$  in general.

\begin{exa}\label{exa-2.1}
Define $\Omega=(0,1)$,  $\cal{F}=\cal{B}(\Omega)$ and $P$ be the
Lebesgue measure on $\Omega$. For $n\in\mathbb{N}$, we define a
random variable $X_n$ by
\begin{eqnarray*}
X_n(\omega)=\left\{
\begin{array}{ll}
1, & \mbox{if}\ \   \omega\in (0,\frac{1}{n^2});\\
0, & \mbox{if}\ \ \omega\in [\frac{1}{n^2}, 1).
\end{array}
\right.
\end{eqnarray*}
Then, for any $p>0$, we have
\begin{eqnarray*}
\sum_{n=1}^{\infty}E[|X_n-0|^p]=\sum_{n=1}^{\infty}\int_0^{\frac{1}{n^2}}1^pdP=\sum_{n=1}^{\infty}\frac{1}{n^2}<\infty,
\end{eqnarray*}
i.e., $X_n\stackrel{S\mbox{-}L^p}{\longrightarrow} 0$.  Obviously,
we have $\|X_n-0\|_{\infty}=1$ for any $n\in\mathbb{N}$. Hence
 $\|X_n-0\|_{\infty}\nrightarrow 0$, which implies that $X_n\stackrel{S\mbox{-}L^{\infty}}{\nrightarrow} 0$.
\end{exa}

The following example shows that
$X_n\stackrel{S\mbox{-}L^{\alpha}}{\longrightarrow} X$ is stronger
than $X_n\stackrel{S_{\alpha}\mbox{-}a.s.}{\longrightarrow} X$  in
general.

\begin{exa}\label{exa-2.2}
Let $\alpha>0$. Define $\Omega=(0,1)$,  $\cal{F}=\cal{B}(\Omega)$
and $P$ be the Lebesgue measure on $\Omega$. For $n\in\mathbb{N}$,
we define a random variable $X_n$ by
\begin{eqnarray*}
X_n(\omega)=\left\{
\begin{array}{ll}
1, & \mbox{if}\ \   \omega\in (0,\frac{1}{n});\\
0, & \mbox{if}\ \ \omega\in [\frac{1}{n}, 1).
\end{array}
\right.
\end{eqnarray*}
It is easy to check that
$X_n\stackrel{S_{\alpha}\mbox{-}a.s.}{\longrightarrow} 0$ but
$X^{\alpha}_n\stackrel{c.c.}{\nrightarrow} 0$ and hence
$X_n\stackrel{S\mbox{-}L^{\alpha}}{\nrightarrow} 0$.
\end{exa}

\begin{rem}
(i) By Examples \ref{exa-2.1} and \ref{exa-2.2} we know that
$X_n\stackrel{S\mbox{-}L^{\infty}}{\longrightarrow} X$ is stronger
than $X_n\stackrel{S_{\alpha}\mbox{-}a.s.}{\longrightarrow} X$  in
general.

(ii) By Example \ref{exa-2.2} and Theorem \ref{pro-2.5}, we know
that $X_n\stackrel{S_{\alpha}\mbox{-}a.s.}{\longrightarrow}
X\nRightarrow X_n\stackrel{S_2\mbox{-}d}{\longrightarrow} X$ in
general.
\end{rem}

Example \ref{exa-2.2}  shows that
$X_n\stackrel{S_{\alpha}\mbox{-}a.s.}{\longrightarrow} X$ does not
imply that $X_n\stackrel{c.c.}{\longrightarrow} X$ in general.
Conversely, the following example   shows that
$X_n\stackrel{c.c.}{\longrightarrow} X$ does not imply
$X_n\stackrel{S_{\alpha}\mbox{-}a.s.}{\longrightarrow} X$ in
general too.

\begin{exa}\label{exa-2.3} Let $\alpha>0$.
Define $\Omega=(0,1)$,  $\cal{F}=\cal{B}(\Omega)$ and $P$ be the
Lebesgue measure on $\Omega$. For $n\in\mathbb{N}$, we define a
random variable $X_n$ by
\begin{eqnarray*}
X_n(\omega)=\left\{
\begin{array}{ll}
1, & \mbox{if}\ \   \omega\in (0,\frac{1}{n^2});\\
\frac{1}{n^{1/\alpha}}, & \mbox{if}\ \ \omega\in [\frac{1}{n^2},
1).
\end{array}
\right.
\end{eqnarray*}

For any $\epsilon>0$, there exists $N$ such that
$\frac{1}{N^{1/\alpha}}<\epsilon$. Then, for any $n\geq N$, we
have $\frac{1}{n^{1/\alpha}}\leq \frac{1}{N^{1/\alpha}}<\epsilon$
and thus
$$
\sum_{n=1}^{\infty}P\{|X_n-0|\geq \epsilon\}\leq \sum_{n=1}^{N-1}P\{|X_n-0|\geq \epsilon\}+\sum_{n=N}^{\infty}\frac{1}{n^2}<\infty.
$$
Hence $X_n\stackrel{c.c.}{\longrightarrow} 0$.

Obviously, for any $\omega\in (0,1)$, we have
$\sum_{n=1}^{\infty}|X_n-0|^{\alpha}=\sum_{n=1}^{\infty}X_n^{\alpha}=\infty$.
Thus $X_n\stackrel{S_{\alpha}\mbox{-}a.s.}{\nrightarrow} 0$.
\end{exa}

The following example shows that if $X$ is nondegenerate, then for
$i\in \{1,2\}$, we do not have $X_n\stackrel{S_i\mbox{-}d}{\longrightarrow} X\Rightarrow X_n\stackrel{P}{\longrightarrow} X.$

\begin{exa}\label{exa-2.4}
Let $\{X,X_n,n\geq 1\}$ be a sequence of i.i.d. random variables.
Obviously, we have $X_n\stackrel{S_i\mbox{-}d}{\longrightarrow} X$
for  $i\in \{1,2\}$. Suppose that $X$ is nondegenerate. Then,
there exists positive constants $c_1$ and $c_2$ such that
$P\{|X_n-X|\geq c_1\}\equiv c_2$ for any $n\in\mathbb{N}$.
Therefore, $X_n\stackrel{P}{\nrightarrow} X$.
\end{exa}

\section{Some open questions}

In this section, we present some open questions for further
research.

\noindent {\bf Question 1.} What is the relation between the
$S_1\mbox{-}d$ convergence and the $S_2\mbox{-}d$ convergence?

\noindent {\bf Question 2.} Does
$X_n\stackrel{S\mbox{-}L^{\infty}}{\longrightarrow} X$ imply that
$X_n\stackrel{S_2\mbox{-}d}{\longrightarrow} X$?

\noindent {\bf Question 3.} Does
$X_n\stackrel{S\mbox{-}L^1}{\longrightarrow} X$ imply that
$X_n\stackrel{S_2\mbox{-}d}{\longrightarrow} X$?

\noindent {\bf Question 4.} Does
$X_n\stackrel{S_{\alpha}\mbox{-}a.s.}{\longrightarrow} X$
($\alpha>0$) imply that
$X_n\stackrel{S_1\mbox{-}d}{\longrightarrow} X$?
 % (By Remark 2.10(ii), we know that $X_n\stackrel{S_{\alpha}\mbox{-}a.s.}{\longrightarrow}
%X\nRightarrow X_n\stackrel{S_2\mbox{-}d}{\longrightarrow} X$ in
%general.)

\noindent {\bf Question 5.} Does
$X_n\stackrel{c.c.}{\longrightarrow} X$
($\alpha>0$) imply that
$X_n\stackrel{S_i\mbox{-}d}{\longrightarrow} X$ for $i\in \{1,2\}$?

\noindent {\bf Question 6.} Can we give a Skorokhod-type theorem
for the strong convergence in distribution and the
$S_{\alpha}\mbox{-}a.s.$ convergence?

\bigskip

\bigskip

{ \noindent {\bf\large Acknowledgments} \quad   This work was
supported by National Natural Science Foundation of China (Grant
No. 11771309), Natural Sciences and Engineering Research
Council of Canada (Grant No. 311945-2013) and the Fundamental Research Funds for the Central Univesities of China.}

\end{document}